\documentclass[secthm,seceqn,amsthm,ussrhead,12pt,reqno]{amsart}
\usepackage[utf8]{inputenc}
\usepackage[english]{babel}
\usepackage{amssymb,amsmath,amsthm,amsfonts,xcolor,enumerate,hyperref,comment,longtable,cleveref}

\usepackage{times}
\usepackage{cite}
\usepackage{pdflscape}
\usepackage{ulem}
\usepackage[mathcal]{euscript}
\usepackage{tikz}
\usepackage{hyperref}
\usepackage{cancel}
\usepackage{stmaryrd}

\usetikzlibrary{arrows}

\newtheorem{theorem}{Theorem}
\newtheorem{lemma}[theorem]{Lemma}

\newtheorem{remark}[theorem]{Remark}

\newtheorem{definition}[theorem]{Definition}

\crefrangeformat{equation}{#3(#1)#4--#5(#2)#6}

\crefname{theorem}{Theorem}{Theorems}
\crefname{lemma}{Lemma}{Lemmas}
\crefname{corollary}{Corollary}{Corollaries}
\crefname{proposition}{Proposition}{Propositions}
\crefname{definition}{Definition}{Definitions}
\crefname{example}{Example}{Examples}
\crefname{remark}{Remark}{Remarks}
\crefname{section}{Section}{Sections}
\crefname{equation}{\unskip}{\unskip}
\crefname{enumi}{\unskip}{\unskip}

\newenvironment{Proof}[1][Proof.]{\begin{trivlist}
\item[\hskip \labelsep {\bfseries #1}]}{\flushright
$\Box$\end{trivlist}}

\usepackage{stmaryrd}
\usepackage{xcolor}

\newcommand{\af}{\alpha}

\begin{document}
\noindent{\Large 
The algebraic classification of nilpotent Tortkara algebras}\footnote{
The work was supported by RSF 18-71-10007. 
The authors thank  Pilar P\'{a}ez-Guill\'{a}n for some constructive comments.
}$^,$\footnote{Corresponding Author: kaygorodov.ivan@gmail.com  }

   \

   {\bf  Ilya Gorshkov$^{a,b,c}$,  Ivan Kaygorodov$^{c}$ \&  Mykola Khrypchenko$^{d}$ \\

    \medskip
}

{\tiny

$^{a}$ Sobolev Institute of Mathematics, Novosibirsk, Russia

$^{b}$ Siberian Federal University, Krasnoyarks, Russia

$^{c}$ CMCC, Universidade Federal do ABC. Santo Andr\'e, Brazil

$^{d}$ Departamento de Matem\'atica, Universidade Federal de Santa Catarina, Florian\'opolis, Brazil
\

\smallskip

   E-mail addresses:

\smallskip
    Ilya Gorshkov (ilygor8@gmail.com)
    
    Ivan Kaygorodov (kaygorodov.ivan@gmail.com) 
    
    Mykola Khrypchenko (nskhripchenko@gmail.com)

}

\ 

\

\ 

\noindent {\bf Abstract.}
{\it We classify all $6$-dimensional nilpotent Tortkara algebras over $\mathbb C$.}

\ 


\ 


\section*{Introduction}

Algebraic classification (up to isomorphism) of $n$-dimensional algebras from a certain variety
defined by some family of polynomial identities is a classical problem in the theory of non-associative algebras.
There are many results related to algebraic classification of small dimensional algebras in varieties of 
Jordan, Lie, Leibniz, Zinbiel and many another algebras \cite{ack, cfk19, gkks, degr3, usefi1, degr1, degr2, ha16, hac18, kv16}.
Another interesting direction in classifications of algebras is geometric classification. We refer the reader to \cite{kppv, kpv, kv16, S90} for results in this direction for 
Jordan, Lie, Leibniz, Zinbiel and other algebras.
In the present paper, we give algebraic classification of nilpotent algebras of a new class of non-associative algebras introduced by Dzhumadildaev in \cite{dzhuma}.

An algebra $\bf A$ is called a {\it Zinbiel algebra} if it satisfies the identity 
$$(xy)z=x(yz+zy).$$ 
Zinbiel algebras were introduced by Loday in \cite{loday} and studied in \cite{cam13, dok, dzhuma5, dzhuma19, kppv, ualbay, yau}. Under the Koszul duality, the operad of Zinbiel algebras is dual to the operad of Leibniz algebras. Zinbiel algebras are also related to Tortkara algebras \cite{dzhuma} and Tortkara triple systems \cite{brem}. More precisely, every Zinbiel algebra with the commutator multiplication gives a Tortkara algebra. An anticommutative algebra $\bf A$ is called a {\it Tortkara algebra} if it satisfies the identity 
$$(ab)(cb)=J(a,b,c)b, \mbox{ where } J(a,b,c)=(ab)c+(bc)a+(ca)b.$$
It is easy to see that every metabelian Lie algebra (i.e., $(xy)(zt)=0$) is a Tortkara algebra
and every dual Mock-Lie algebra (i.e., antiassociative and anticommutative) is a Tortkara algebra.

Our method of classification of nilpotent Tortkara algebras is based on calculations of central extensions of smaller nilpotent algebras from the same variety.
Central extensions play an important role in quantum mechanics: one of the earlier
encounters is by means of Wigner's theorem which states that a symmetry of a quantum
mechanical system determines an (anti-)unitary transformation of a Hilbert space.
Another area of physics where one encounters central extensions is the quantum theory
of conserved currents of a Lagrangian. These currents span an algebra which is closely
related to so called affine Kac-Moody algebras, which are universal central extensions
of loop algebras.
Central extensions are needed in physics, because the symmetry group of a quantized
system usually is a central extension of the classical symmetry group, and in the same way
the corresponding symmetry Lie algebra of the quantum system is, in general, a central
extension of the classical symmetry algebra. Kac-Moody algebras have been conjectured
to be a symmetry groups of a unified superstring theory. The centrally extended Lie
algebras play a dominant role in quantum field theory, particularly in conformal field
theory, string theory and in $M$-theory.
In the theory of Lie groups, Lie algebras and their representations, a Lie algebra extension
is an enlargement of a given Lie algebra $g$ by another Lie algebra $h.$ Extensions
arise in several ways. There is a trivial extension obtained by taking a direct sum of
two Lie algebras. Other types are a split extension and a central extension. Extensions may
arise naturally, for instance, when forming a Lie algebra from projective group representations.
A central extension and an extension by a derivation of a polynomial loop algebra
over finite-dimensional simple Lie algebra give a Lie algebra which is isomorphic to a
non-twisted affine Kac-Moody algebra \cite[Chapter 19]{bkk}. Using the centrally extended loop
algebra one may construct a current algebra in two spacetime dimensions. The Virasoro
algebra is the universal central extension of the Witt algebra, the Heisenberg algebra is
the central extension of a commutative Lie algebra  \cite[Chapter 18]{bkk}.

The algebraic study of central extensions of Lie and non-Lie algebras has a very long history \cite{omirov,ha17,hac16,kkl18,ss78,zusmanovich}.
Thus, Skjelbred and Sund used central extensions of Lie algebras for a classification of nilpotent Lie algebras  \cite{ss78}.
After that, the method introduced by Skjelbred and Sund was used to 
describe  all non-Lie central extensions of  all $4$-dimensional Malcev algebras \cite{hac16},
all non-associative central extensions of $3$-dimensional Jordan algebras \cite{ha17},
all anticommutative central extensions of $3$-dimensional anticommutative algebras \cite{cfk182},
all central extensions   of $2$-dimensional algebras \cite{cfk18}.
The method of central extensions  was used to describe
all $4$-dimensional nilpotent associative algebras \cite{degr1},
all $4$-dimensional nilpotent bicommutative algebras \cite{kpv19},
all $4$-dimensional nilpotent Novikov  algebras \cite{kkk18},
all $5$-dimensional nilpotent Jordan algebras \cite{ha16},
all $5$-dimensional nilpotent restricted Lie algebras \cite{usefi1},
all $6$-dimensional nilpotent Lie algebras \cite{degr3,degr2},
all $6$-dimensional nilpotent Malcev algebras \cite{hac18},
all $6$-dimensional nilpotent binary Lie algebras\cite{ack},  
all $6$-dimensional nilpotent anticommutative $\mathfrak{CD}$-algebras \cite{ack}
and some other.

\section{Preliminaries}
\subsection{Method of classification of nilpotent algebras}
Throughout this paper, we  use the notations and methods described in \cite{ha17,hac16,cfk18}
and adapted for the Tortkara case with some modifications. 
From now on, we will give only some important definitions.

Let $({\bf A}, \cdot)$ be a Tortkara algebra over $\mathbb C$ and $\mathbb V$ a vector space over the same base field. Then the $\mathbb C$-linear space ${\rm Z_T^{2}}\left(
\bf A,\mathbb V \right) $ is defined as the set of all skew-symmetric bilinear maps $%
\theta :{\bf A} \times {\bf A} \longrightarrow {\mathbb V},$
such that \[\theta(xy,zt)+\theta(xt,zy)=\theta(J(x,y,z),t)+ \theta(J(x,t,z),y).\]
 Its elements will be called \textit{cocycles}. For a
linear map $f$ from $\bf A$ to  $\mathbb V$, if we define $\delta f\colon {\bf A} \times
{\bf A} \longrightarrow {\mathbb V}$ by $\delta f \left( x,y \right) =f(xy )$, then $%
\delta f\in {\rm Z_T^{2}}\left( {\bf A},{\mathbb V} \right) $. Let ${\rm B}^{2}\left(
{\bf A},{\mathbb V}\right) =\left\{ \theta =\delta f\ :f\in {\rm Hom}\left( {\bf A},{\mathbb V}\right) \right\} $.
One can easily check that ${\rm B}^{2}(\bf A,\mathbb V)$ is a linear subspace of $%
{\rm Z_T^{2}}\left( {\bf A},{\mathbb V}\right) $ whose elements are called \textit{%
coboundaries}. We define the \textit{second cohomology space} $%
{\rm H_T^{2}}\left( {\bf A},{\mathbb V}\right) $ as the quotient space ${\rm Z_T^{2}}%
\left( {\bf A},{\mathbb V}\right) \big/{\rm B}^{2}\left( {\bf A},{\mathbb V}\right) $.

\bigskip

Let $\operatorname{Aut}\left( {\bf A} \right) $ be the automorphism group of the Tortkara  algebra $%
{\bf A} $ and let $\phi \in \operatorname{Aut} \left( {\bf A}\right) $. For $\theta \in
{\rm Z_T^{2}}\left( {\bf A},{\mathbb V}\right) $ define $\phi \theta \left( x,y\right)
=\theta \left( \phi \left( x\right) ,\phi \left( y\right) \right) $. Then $%
\phi \theta \in {\rm Z_T^{2}}\left( {\bf A},{\mathbb V}\right) $. So, $\operatorname{Aut}\left( {\bf A}\right) $
acts on ${\rm Z_T^{2}}\left( {\bf A},{\mathbb V}\right) $. It is easy to verify that $%
{\rm B}^{2}\left( {\bf A},{\mathbb V}\right) $ is invariant under the action of $\operatorname{Aut}\left(
{\bf A}\right) $, and thus $\operatorname{Aut}\left( {\bf A}\right) $ acts on $%
{\rm H_T^{2}}\left( {\bf A},{\mathbb V}\right) $.

\bigskip

Let $\bf A$ be a Tortkara  algebra of dimension $m<n$ over $\mathbb C$, and ${\mathbb V}$ be a $\mathbb C$-vector
space of dimension $n-m$. For any $\theta \in {\rm Z_T^{2}}\left(
{\bf A},{\mathbb V}\right) $ define on the linear space ${\bf A}_{\theta }:={\bf A}\oplus {\mathbb V}$ the
bilinear product ``$\left[ -,-\right] _{{\bf A}_{\theta }}$'' by $%
\left[ x+x^{\prime },y+y^{\prime }\right] _{{\bf A}_{\theta }}=
 xy +\theta \left( x,y\right) $ for all $x,y\in {\bf A},x^{\prime },y^{\prime }\in {\mathbb V}$.
Then ${\bf A}_{\theta }$ is a Tortkara algebra called an $(n-m)$%
{\it{-dimensional central extension}} of ${\bf A}$ by ${\mathbb V}$. In fact, ${\bf A_{\theta}}$ is a Tortkara algebra \textit{if and only if} $\theta \in {\rm Z_T}^2({\bf A}, {\mathbb V})$.

We also call the
set $\operatorname{Ann}(\theta)=\left\{ x\in {\bf A}:\theta \left( x, {\bf A} \right) =0\right\} $
the {\it{annihilator}} of $\theta $. We recall that the {\it{annihilator}} of an  algebra ${\bf A}$ is defined as
the ideal $\operatorname{Ann}\left( {\bf A} \right) =\left\{ x\in {\bf A}:  x{\bf A}=0\right\}$ and observe
 that
$\operatorname{Ann}\left( {\bf A}_{\theta }\right) = \operatorname{Ann}(\theta) \cap \operatorname{Ann}\left( {\bf A}\right)
 \oplus {\mathbb V}.$

\medskip

We have the next  key result:

\begin{lemma}
Let ${\bf A}$ be an $n$-dimensional Tortkara algebra such that $\dim(\operatorname{Ann}({\bf A}))=m\neq0$. Then there exists, up to isomorphism, a unique $(n-m)$-dimensional Tortkara algebra ${\bf A}^{\prime }$ and a bilinear map $\theta \in {\rm Z_T}^2({\bf A}, {\mathbb V})$ with $\operatorname{Ann}({\bf A})\cap \operatorname{Ann}(\theta)=0$, where $\mathbb V$ is a vector space of dimension m, such that ${\bf A}\cong {\bf A}^{\prime }_{\theta}$ and $
{\bf A}/\operatorname{Ann}({\bf A})\cong {\bf A}^{\prime }$.
\end{lemma}

\begin{Proof}

Let ${\bf A}^{\prime }$ be a linear complement of $\operatorname{Ann}({\bf A})$ in ${\bf A}$. Define a linear map $P: {\bf A} \longrightarrow {\bf A}^{\prime }$ by $P(x+v)=x$ for $x\in {\bf A}^{\prime }$ and $v\in \operatorname{Ann}({\bf A})$ and define a multiplication on ${\bf A}^{\prime }$ by $[x, y]_{{\bf A}^{\prime }}=P(x y)$ for $x, y \in {\bf A}^{\prime }$.
For $x, y \in {\bf A}$ we have
$$P(xy)=P((x-P(x)+P(x))(y- P(y)+P(y)))=P(P(x) P(y))=[P(x), P(y)]_{{\bf A}^{\prime }}$$

Since $P$ is a homomorphism, then $P({\bf A})={\bf A}^{\prime }$ is a Tortkara algebra and $
{\bf A}/\operatorname{Ann}({\bf A})\cong {\bf A}^{\prime },$ which gives us the uniqueness of ${\bf A}^{\prime }$. Now, define the map $\theta: {\bf A}^{\prime }\times{\bf A}^{\prime }\longrightarrow \operatorname{Ann}({\bf A})$ by $\theta(x,y)=xy-[x,y]_{{\bf A}^{\prime }}$. Then ${\bf A}^{\prime }_{\theta}$ is ${\bf A}$ and therefore $\theta \in {\rm Z_T}^2({\bf A}, {\mathbb V})$ and $\operatorname{Ann}({\bf A})\cap \operatorname{Ann}(\theta)=0$.
\end{Proof}

\bigskip

However, in order to solve the isomorphism problem, we need to study the
action of $\operatorname{Aut}\left( {\bf A}\right) $ on ${\rm H_T^{2}}\left( {\bf A},{\mathbb V}%
\right) $. To this end, let us fix $e_{1},\ldots ,e_{s}$ a basis of ${\mathbb V}$, and $%
\theta \in {\rm Z_T^{2}}\left( {\bf A},{\mathbb V}\right) $. Then $\theta $ can be uniquely
written as $\theta \left( x,y\right) =\underset{i=1}{\overset{s}{\sum }}%
\theta _{i}\left( x,y\right) e_{i}$, where $\theta _{i}\in
{\rm Z_T^{2}}\left( {\bf A},\mathbb C\right) $. Moreover, $\operatorname{Ann}(\theta)=\operatorname{Ann}(\theta _{1})\cap \operatorname{Ann}(\theta _{2})\cdots \cap \operatorname{Ann}(\theta _{s})$. Further, $\theta \in
{\rm B}^{2}\left( {\bf A},{\mathbb V}\right) $\ if and only if all $\theta _{i}\in {\rm B}^{2}\left( {\bf A},%
\mathbb C\right) $.

\bigskip

\begin{definition}
Given a Tortkara algebra ${\bf A}$, if ${\bf A}=I\oplus \mathbb Cx$
is a direct sum of two ideals, then $\mathbb Cx$ is called an {\it{%
annihilator component}} of ${\bf A}$. 
\end{definition}

\begin{definition}
A central extension of an algebra $\bf A$ without an annihilator component is called a non-split central extension.
\end{definition}

It is not difficult to prove (see \cite[%
Lemma 13]{hac16}) that, given a Tortkara algebra ${\bf A}_{\theta}$, if we write as
above $\theta \left( x,y\right) =\underset{i=1}{\overset{s}{\sum }}$
 $\theta_{i}\left( x,y\right) e_{i}\in {\rm Z_T^{2}}\left( {\bf A},{\mathbb V}\right) $ and we have
$\operatorname{Ann}(\theta)\cap \operatorname{Ann}\left( {\bf A}\right) =0$, then ${\bf A}_{\theta }$ has an
annihilator component if and only if $\left[ \theta _{1}\right] ,\left[
\theta _{2}\right] ,\ldots ,\left[ \theta _{s}\right] $ are linearly
dependent in ${\rm H_T^{2}}\left( {\bf A},\mathbb C\right) $.

\bigskip

Let ${\mathbb V}$ be a finite-dimensional vector space over $\mathbb C$. The {\it{%
Grassmannian}} $G_{k}\left( {\mathbb V}\right) $ is the set of all $k$-dimensional
linear subspaces of $ {\mathbb V}$. Let $G_{s}\left( {\rm H_T^{2}}\left( {\bf A},\mathbb C%
\right) \right) $ be the Grassmannian of subspaces of dimension $s$ in $%
{\rm H_T^{2}}\left( {\bf A},\mathbb C\right) $. There is a natural action of $%
\operatorname{Aut} \left( {\bf A}\right) $ on $G_{s}\left( {\rm H_T^{2}}\left( {\bf A},\mathbb C%
\right) \right) $. Let $\phi \in \operatorname{Aut} \left( {\bf A}\right) $. For $W=\left\langle %
\left[ \theta _{1}\right] ,\left[ \theta _{2}\right] ,\dots,\left[ \theta _{s}%
\right] \right\rangle \in G_{s}\left( {\rm H_T^{2}}\left( {\bf A},\mathbb C%
\right) \right) $ define $\phi W=\left\langle \left[ \phi \theta _{1}\right]
,\left[ \phi \theta _{2}\right] ,\dots,\left[ \phi \theta _{s}\right]
\right\rangle $. Then $\phi W\in G_{s}\left( {\rm H_T^{2}}\left( {\bf A},\mathbb C \right) \right) $. We denote the orbit of $W\in G_{s}\left(
{\rm H_T^{2}}\left( {\bf A},\mathbb C\right) \right) $ under the action of $%
\operatorname{Aut} \left( {\bf A}\right) $ by $\mathrm{Orb}\left( W\right) $. Since, given
\begin{equation*}
W_{1}=\left\langle \left[ \theta _{1}\right] ,\left[ \theta _{2}\right] ,\dots,%
\left[ \theta _{s}\right] \right\rangle ,W_{2}=\left\langle \left[ \vartheta
_{1}\right] ,\left[ \vartheta _{2}\right] ,\dots,\left[ \vartheta _{s}\right]
\right\rangle \in G_{s}\left( {\rm H_T^{2}}\left( {\bf A},\mathbb C\right)\right),
\end{equation*}%
we easily have that if $W_{1}=W_{2}$, then $\underset{i=1}{\overset{s}{%
\cap }}\operatorname{Ann}(\theta _{i})\cap \operatorname{Ann}\left( {\bf A}\right) =\underset{i=1}{\overset{s}%
{\cap }}\operatorname{Ann}(\vartheta _{i})\cap \operatorname{Ann}\left( {\bf A}\right) $, we can introduce
the set

\begin{equation*}
{\bf T}_{s}\left( {\bf A}\right) =\left\{ W=\left\langle \left[ \theta _{1}\right] ,%
\left[ \theta _{2}\right] ,\dots,\left[ \theta _{s}\right] \right\rangle \in
G_{s}\left( {\rm H_T^{2}}\left( {\bf A},\mathbb C\right) \right) :\underset{i=1}{%
\overset{s}{\cap }}\operatorname{Ann}(\theta _{i})\cap \operatorname{Ann}\left( {\bf A}\right) =0\right\},
\end{equation*}
which is stable under the action of $\operatorname{Aut}\left( {\bf A}\right) $.

\medskip

Now, let ${\mathbb V}$ be an $s$-dimensional linear space and let us denote by $%
{\bf E}\left( {\bf A},{\mathbb V}\right) $ the set of all {\it non-split $s$-dimensional} central extensions of ${\bf A}$ by
${\mathbb V}.$ We can write
\begin{equation*}
{\bf E}\left( {\bf A},{\mathbb V}\right) =\left\{ {\bf A}_{\theta }:\theta \left( x,y\right) =\underset{%
i=1}{\overset{s}{\sum }}\theta _{i}\left( x,y\right) e_{i}\mbox{
and }\left\langle \left[ \theta _{1}\right] ,\left[ \theta _{2}\right] ,\dots,%
\left[ \theta _{s}\right] \right\rangle \in {\bf T}_{s}\left( {\bf A}\right) \right\} .
\end{equation*}%
Also, we have the next result, which can be proved as \cite[Lemma 17]{hac16}.

\begin{lemma}
 Let ${\bf A}_{\theta },{\bf A}_{\vartheta }\in {\bf E}\left( {\bf A},{\mathbb V}\right).$ 
 Suppose that $\theta \left( x,y\right) =\underset{i=1}{\overset{s}{\sum }}%
\theta _{i}\left( x,y\right) e_{i}$ and $\vartheta \left( x,y\right) =%
\underset{i=1}{\overset{s}{\sum }}\vartheta _{i}\left( x,y\right) e_{i}$.
Then the Tortkara algebras ${\bf A}_{\theta }$ and ${\bf A}_{\vartheta } $ are isomorphic
if and only if 
\[\mathrm{Orb}\left\langle \left[ \theta _{1}\right] ,%
\left[ \theta _{2}\right] ,\dots,\left[ \theta _{s}\right] \right\rangle =%
\mathrm{Orb}\left\langle \left[ \vartheta _{1}\right] ,\left[ \vartheta
_{2}\right] ,\dots,\left[ \vartheta _{s}\right] \right\rangle. \]
\end{lemma}

Thus, there exists a one-to-one correspondence between the set of $\operatorname{Aut}
\left( {\bf A}\right) $-orbits on ${\bf T}_{s}\left( {\bf A}\right) $ and the set of
isomorphism classes of ${\bf E}\left( {\bf A},{\mathbb V}\right) $. Consequently, we have a
procedure that allows us, given a Tortkara algebra ${\bf A}^{\prime }$ of
dimension $n$, to construct all non-split central extensions of ${\bf A}^{\prime }.$ This procedure is as follows:

\medskip

{\centerline{\it Procedure}}

\begin{enumerate}
\item For a given (nilpotent) Tortkara algebra $\bf{A}^{\prime }$
of dimension $n-s$, determine ${\bf T}_{s}(\bf{A}^{\prime })$
and $\operatorname{Aut}(\bf{A}^{\prime })$.

\item Determine the set of $\operatorname{Aut}(\bf{A}^{\prime })$-orbits on $%
{\bf T}_{s}(\bf{A}^{\prime })$.

\item For each orbit, construct the Tortkara algebra corresponding to one of its
representatives.
\end{enumerate}

The above described method gives all (Malcev and non-Malcev) Tortkara 
algebras. But we are interested in developing this method in such a way
that it only gives non-Malcev Tortkara algebras, because the classification of all Malcev-Tortkara algebras is a part of  the classification of all Malcev algebras \cite{kpv, hac18}.
Clearly, any central extension of a non-Malcev Tortkara  algebra is non-Malcev. 
But a Malcev-Tortkara algebra may have extensions which are Malcev algebras. More precisely, 
let ${\bf M}$ be
a Malcev algebra and $\theta \in {\rm Z_T^{2}}\left( {\bf M},%
{\mathbb C}\right) $. Then ${\bf M}_{\theta }$ is a Malcev algebra if and
only if 
\begin{equation*}
\theta \left(  wy, xz \right) =
\theta \left( (wx)y, z \right)+
\theta \left( (xy)z, w \right)+
\theta \left( (yz)w, x \right)+
\theta \left( (zw)x, y \right) 
\end{equation*}%
for all $x,y,z,w\in {\bf M}$. Define the subspace 
${\rm Z_{TM}^{2}}\left( {\bf M},{\mathbb C}\right) $ of ${\rm Z_{T}^{2}}\left( {\bf M},{\mathbb C}\right) $ by%
\begin{equation*}
{\rm Z_{TM}^{2}}\left( {\bf M},{\mathbb C}\right) =\left\{ 
\begin{array}{c}
\theta \in {\rm Z_{T}^{2}}\left( {\bf M},{\mathbb C}\right) :
\theta \left(  wy, xz \right) =
\theta \left( (wx)y, z \right)+
\theta \left( (xy)z, w \right)+\\
\theta \left( (yz)w, x \right)+
\theta \left( (zw)x, y \right), \text{ for all } x, y,z,w\in {\bf M}%
\end{array}%
\right\} .
\end{equation*}
Observe that ${\rm B}^{2}\left( {\bf M},{\mathbb C}\right)\subseteq{\rm Z_{TM}^{2}}\left( {\bf M},{\mathbb C}\right)$.
Let ${\rm H_{TM}^{2}}\left( {\bf M},{\mathbb C}\right) =
{\rm Z_{TM}^{2}}\left( {\bf M},{\mathbb C}\right) \big/B%
^{2}\left( {\bf M},{\mathbb C}\right) $. Then ${\rm H_{TM}^{2}}\left( {\bf M},{\mathbb C}\right) $ is a subspace of $%
{\rm H_{T}^{2}}\left( {\bf M},{\mathbb C}\right) $. Define 
\begin{eqnarray*}
{\bf R}_{s}\left( {\bf M}\right)  &=&\left\{ {\bf W}\in {\bf 
T}_{s}\left( {\bf M}\right) :{\bf W}\in G_{s}\left( {\rm H_{TM}^{2}}\left( {\bf M},{\mathbb C}\right) \right) \right\} , \\
{\bf U}_{s}\left( {\bf M}\right)  &=&\left\{ {\bf W}\in {\bf
T}_{s}\left( {\bf M}\right) :{\bf W}\notin G_{s}\left( {\rm H
_{TM}^{2}}\left( {\bf M},{\mathbb C}\right) \right) \right\} .
\end{eqnarray*}%
Then ${\bf T}_{s}\left( {\bf M}\right) ={\bf R}_{s}\left( 
{\bf M}\right) $ $\mathbin{\mathaccent\cdot\cup}$ ${\bf U}_{s}\left( 
{\bf M}\right) $. The sets ${\bf R}_{s}\left( {\bf M}\right) $
and ${\bf U}_{s}\left( {\bf M}\right) $ are stable under the action
of $\operatorname{Aut}\left( {\bf M}\right) $. Thus, the Tortkara  algebras
corresponding to the representatives of $\operatorname{Aut}\left( {\bf M}\right) $%
-orbits on ${\bf R}_{s}\left( {\bf M}\right) $ are Malcev algebras,
while those corresponding to the representatives of $\operatorname{Aut}\left( {\bf M}%
\right) $-orbits on ${\bf U}_{s}\left( {\bf M}\right) $ are not
Malcev algebras. Hence, we may construct all non-split non-Malcev Tortkara algebras $%
\bf{A}$ of dimension $n$ with $s$-dimensional annihilator 
from a given Tortkara algebra $\bf{A}%
^{\prime }$ of dimension $n-s$ in the following way:

\begin{enumerate}
\item If $\bf{A}^{\prime }$ is non-Malcev, then apply the Procedure.

\item Otherwise, do the following:

\begin{enumerate}
\item Determine ${\bf U}_{s}\left( \bf{A}^{\prime }\right) $ and $%
\operatorname{Aut}(\bf{A}^{\prime })$.

\item Determine the set of $\operatorname{Aut}(\bf{A}^{\prime })$-orbits on ${\bf U%
}_{s}\left( \bf{A}^{\prime }\right) $.

\item For each orbit, construct the Tortkara  algebra corresponding to one of its
representatives.
\end{enumerate}
\end{enumerate}

\medskip

\subsection{Notations}
Let us introduce the following notations. Let ${\bf A}$ be a Tortkara algebra with
a basis $e_{1},e_{2}, \ldots, e_{n}$. Then by $\Delta _{ij}$\ we will denote the
bilinear form
$\Delta _{ij}:{\bf A}\times {\bf A}\longrightarrow \mathbb C$
with $\Delta _{ij}\left( e_{l},e_{m}\right) =0$, if $%
\left\{ i,j\right\} \neq \left\{ l,m\right\} $,
and $\Delta _{ij}\left( e_{i},e_{j}\right)=-\Delta _{ij}\left( e_{j},e_{i}\right)=1$. The set $\left\{ \Delta
_{ij}:1\leq i < j\leq n\right\} $ is a basis for the linear space of skew-symmetric
bilinear forms on ${\bf A}$, so every $\theta \in
{\rm Z_T^{2}}\left( {\bf A},\bf \mathbb V \right) $ can be uniquely written as $%
\theta =\underset{1\leq i<j\leq n}{\sum }c_{ij}\Delta _{ij}$, where $%
c_{ij}\in \mathbb C$. We also denote by

$$\begin{array}{lll}
{\mathbb T}^i_j& \mbox{the }j\mbox{th }i\mbox{-dimensional nilpotent Tortkara algebra}, \\
{\mathfrak{N}}_i& \mbox{the }i\mbox{-dimensional algebra with zero product}. \\
\end{array}$$

\subsection{Central extensions of nilpotent low dimensional Tortkara algebras}
There are no nontrivial $1$- and $2$-dimensional nilpotent Tortkara algebras.
Thanks to \cite{cfk182} we have the description of all non-split $3$- and $4$-dimensional  nilpotent Tortkara algebras:

$\begin{array}{lllll lll}
{\mathbb T}^3_{01}&:& 
e_1e_2=e_3;\\
{\mathbb T}^4_{02} &:& 
e_1e_2=e_3, & e_1e_3=e_4. 
\end{array}$

Also, it is easy to see that every anticommutative central extension of $\mathfrak{N}_i$ is a Lie algebra.

\section{Central extensions of nilpotent $4$-dimensional Tortkara algebras}

\subsection{The algebraic classification of $4$-dimensional nilpotent Tortkara algebras}

\begin{equation*}
\begin{array}{|l|l|l|l|} 
\hline
\mbox{$\bf A$}  & \mbox{ multiplication table } & \mbox{${\rm {\rm H_{TM}^2}}({\bf A})$} & \mbox{${\rm {\rm H_{T}^2}}({\bf A})$} \\ 

\hline
\hline

{\mathbb T}^4_{01}& e_1e_2 = e_3 &
\langle [\Delta_{13}], [\Delta_{14}], [\Delta_{23}], [\Delta_{24}], [\Delta_{34}]   \rangle  &  \mbox{${\rm {\rm H_{TM}^2}}({{\mathbb T}^4_{01}})$}  \\
\hline
{\mathbb T}^4_{02}& e_1e_2=e_3, e_1e_3=e_4 &
\langle [\Delta_{14}], [\Delta_{23}] \rangle & 
\mbox{${\rm {\rm H_{TM}^2}}({\mathbb T}^4_{01})$} \oplus \langle  [\Delta_{24}]\rangle  \\
\hline

\end{array}
\end{equation*}





\subsection{Central extensions of ${\mathbb T}^4_{02}$}
Since ${\rm H_{TM}^2}({{\mathbb T}^4_{02}}) \neq {\rm H_{T}^2}({{\mathbb T}^4_{02}}),$
 we can find some Tortkara (non-Malcev) algebras.
 Let us introduce the following notations
 \[  \nabla_1= [\Delta_{14}], \ \nabla_2=[\Delta_{23}], \ \nabla_3=[\Delta_{24}].\] 
The automorphism group of ${\mathbb T}^4_{02}$ 
consists of the invertible
matrices of the form 
\[\phi = \begin{pmatrix} 
x& 0 & 0 & 0\\
z & y & 0 & 0 \\
u & v & xy  & 0\\
h & g & xv & x^2y
\end{pmatrix}.\]
 Since 
\[
\phi^t
\begin{pmatrix} 
0 & 0 & 0 & \alpha_1\\
0 & 0 & \alpha_2 & \alpha_3 \\
0 & -\alpha_2 & 0  & 0\\
-\alpha_1  & -\alpha_3 & 0& 0
\end{pmatrix}
\phi=
\begin{pmatrix}
0& \alpha^*& \alpha^{**}&\alpha_1^*\\
-\alpha^*& 0& \alpha_2^*&\alpha_3^*\\
-\alpha^{**}& -\alpha_2^*& 0& 0\\
-\alpha_1^*& -\alpha_3^*& 0& 0

\end{pmatrix},
\]
where
$$\begin{array}{lcl}
\alpha_1^*&=&(\alpha_1x+\alpha_3z)x^2y,\\
\alpha_2^*&=&(\alpha_2y+\alpha_3v)xy,\\
\alpha_3^*&=&\alpha_3x^2y^2,
\end{array}$$
we obtain that the action of $\operatorname{Aut}\left(\mathbb T_{02}^4 \right)$ on the subspace $\Big\langle \sum\limits_{i=1}^3 \alpha_i \nabla_i\Big\rangle$ is given by
$\Big\langle \sum\limits_{i=1}^3 \alpha_i^* \nabla_i\Big\rangle.$

\subsubsection{$1$-dimensional central extensions of ${\mathbb T}^4_{02}$}
Consider a cocycle
$$
\theta=\alpha_1 \nabla_1+\alpha_2  \nabla_2 +\alpha_3  \nabla_3.
$$ 
To construct non-Malcev Tortkara algebras, we will be interested only in $\theta$ with $\alpha_3\neq0.$
Then  
choosing $x=1, y=\frac{1}{\sqrt{\alpha_3}}, z=-\frac{\alpha_1}{\alpha_3}, v=-\frac{\alpha_2}{\sqrt{(\alpha_3)^3}}$
we have the representative $\big \langle [\Delta_{24}] \big \rangle.$ 
It gives the following new Tortkara algebra:
$$\begin{array}{lllll lll}
{\mathbb T}_{10}^5&:&   e_1e_2=e_3,& e_1e_3=e_4,& e_2e_4=e_5.  \\
\end{array}$$

\subsubsection{$2$-dimensional central extensions of ${\mathbb T}^4_{02}$}

Consider the vector space generated by the following two cocycles 
$$\begin{array}{lcl}
\theta_1&=&\alpha_1  \nabla_1+\alpha_2  \nabla_2+\alpha_3 \nabla_3,\\
\theta_2&=&\beta_1 \nabla_1+\beta_2 \nabla_2.\\
\end{array}$$
We are interested in $\theta_1$ with $\alpha_3\neq 0$. Then
\begin{enumerate}
    \item if $\beta_1\neq0, \beta_2\neq0,$
    then choosing $y=\frac{\beta_1x^2}{\beta_2}, z=-\frac{\alpha_1x}{\alpha_3}, v=-\frac{\alpha_2y}{\alpha_3},$
    we have the representative $\langle \nabla_1+\nabla_2, \ \nabla_3\rangle.$

    \item if $\beta_1=0, \beta_2\neq 0,$ then choosing $z=-\frac{\alpha_1x}{\alpha_3}, v=-\frac{\alpha_2y}{\alpha_3},$
    we have the representative $\langle \nabla_2, \ \nabla_3\rangle.$
    
    \item if $\beta_1\neq 0, \beta_2 = 0,$ then choosing $z=-\frac{\alpha_1x}{\alpha_3}, v=-\frac{\alpha_2y}{\alpha_3},$
    we have the representative $\langle \nabla_1, \ \nabla_3\rangle.$
        
\end{enumerate}

It is easy to see that these three orbits are pairwise distinct.
So, we have all non-split non-Malcev $2$-dimensional central extensions of ${\mathbb T}^4_{02}:$
$$\begin{array}{lllll lll}
{\mathbb T}_{01}^6&:&  
e_1e_2=e_3,& e_1e_3=e_4,& e_1e_4=e_5,& e_2e_3=e_5,& e_2e_4=e_6;  \\

{\mathbb T}_{02}^6&:& 
e_1e_2=e_3,& e_1e_3=e_4,& e_2e_3=e_5, & e_2e_4=e_6;  \\

{\mathbb T}_{03}^6&:&    
e_1e_2=e_3,& e_1e_3=e_4,& e_1e_4=e_5, & e_2e_4=e_6.
\end{array}$$

\section{Central extensions of  $5$-dimensional nilpotent Tortkara algebras}

\subsection{The algebraic classification of $5$-dimensional nilpotent Tortkara algebras}
Since ${\rm H_{T}^2}({{\mathbb T}^4_{01}})={\rm H_{TM}^2}({{\mathbb T}^4_{01}})$, each central extension of ${{\mathbb T}^4_{01}}$ is a Malcev algebra. On the other hand, by some easy calculations, using the algebraic classification of $5$-dimensional nilpotent Malcev algebras \cite{hac18}, one can see that each  $5$-dimensional nilpotent Malcev algebra is a Tortkara algebra. 
Combining this with the result of previous section, we have the algebraic classification of all 
nontrivial $5$-dimensional nilpotent Tortkara algebras:

$$
\begin{array}{|l|l|l|l|} 
\hline
\bf A  & \mbox{ multiplication table } & {\rm H_{TM}^2}({\bf A}) & {\rm H_{T}^2}({\bf A}) \\ 
\hline
\hline
{\mathbb T}^5_{01}& 
\begin{array}{l}
e_1e_2 = e_3 
\end{array}&   

 \Big\langle
\begin{array}{l} 
 [\Delta_{13}] ,  [\Delta_{14}] , [\Delta_{15}] ,  [\Delta_{23}],    [\Delta_{24}], \\ \relax
 [ \Delta_{25} ] ,  [\Delta_{34}], [\Delta_{35}] ,  [\Delta_{45}]
\end{array} 
\Big\rangle

&{\rm H_{TM}^2}({\mathbb T}^5_{01})   \\
\hline
{\mathbb T}^5_{02}& 
\begin{array}{l}
e_1e_2=e_3, e_1e_3=e_4 
\end{array}&
\Big\langle  
[\Delta_{14}], [\Delta_{15}], [\Delta_{23}], [\Delta_{25}],  [\Delta_{35}] 
\Big\rangle
&{\rm H_{TM}^2}({\mathbb T}^5_{02}) \oplus \langle [\Delta_{24}],[\Delta_{45}] \rangle  
\\
\hline

{\mathbb T}^5_{03} & 
\begin{array}{l}
e_1e_2=e_4, e_1e_3=e_5
\end{array}& 
\Big\langle 
\begin{array}{l}
[\Delta_{14}], [\Delta_{15}], [\Delta_{23}], [\Delta_{24}], \\ \relax
[\Delta_{25}], [\Delta_{34}], [\Delta_{35}]
\end{array}
\Big\rangle & {\rm H_{TM}^2}({\mathbb T}^5_{03})\\
\hline

{\mathbb T}^5_{04} & 
\begin{array}{l}
e_1e_2=e_3,  e_1e_3=e_4,\\ e_2e_3=e_5 
\end{array}& 
\Big\langle 
[\Delta_{14}],[\Delta_{15}]+[\Delta_{24}],[\Delta_{25}]
\Big\rangle & {\rm H_{TM}^2}({\mathbb T}^5_{04})\oplus\langle[\Delta_{15}]\rangle\\
\hline

{\mathbb T}^5_{05}& 
\begin{array}{l}
e_1e_2=e_5, e_3e_4=e_5
\end{array}&

\Big\langle 
\begin{array}{l} 
[\Delta_{12}], [\Delta_{13}], [\Delta_{14}], [\Delta_{15}],  [\Delta_{23}],\\ \relax
[\Delta_{24}], [\Delta_{25}], [\Delta_{35}], [\Delta_{45}] 
\end{array} 
\Big\rangle 
&  {\rm H_{TM}^2}({\mathbb T}^5_{05})\\
\hline

{\mathbb T}^5_{06}& 
\begin{array}{l}
e_1e_2=e_3, e_1e_4=e_5,\\
e_2e_3=e_5
\end{array}& 
\Big\langle
[\Delta_{13}],[\Delta_{14}],[\Delta_{24}],[\Delta_{25}],[\Delta_{34}] 
\Big\rangle& {\rm H_{TM}^2}({\mathbb T}^5_{06})\oplus\langle[\Delta_{15}],[\Delta_{45}]\rangle\\
\hline

{\mathbb T}^5_{07}& 
\begin{array}{l}
e_1e_2=e_3,  e_3e_4=e_5 
\end{array}&
\Big\langle 
[\Delta_{13}], [\Delta_{14}], [\Delta_{23}], [\Delta_{24}]
\Big\rangle&  {\rm H_{TM}^2}({\mathbb T}^5_{07}) \\
\hline

{\mathbb T}^5_{08}& 
\begin{array}{l}
e_1e_2=e_3, e_1e_3=e_4, \\ 
e_1e_4=e_5
\end{array}& 
\Big\langle
[\Delta_{15}], [\Delta_{23}]
\Big\rangle & {\rm H_{TM}^2}({\mathbb T}^5_{08})\oplus\langle[\Delta_{24}], [\Delta_{25}]\rangle\\
\hline

{\mathbb T}^5_{09}&
\begin{array}{l}
e_1e_2=e_3, e_1e_3=e_4,\\ 
e_1e_4=e_5,  e_2e_3=e_5
\end{array} & 
\Big\langle
[\Delta_{14}], [\Delta_{15}]+[\Delta_{24}]
\Big\rangle & 
{\rm H_{TM}^2}({\mathbb T}^5_{09})\oplus\langle[\Delta_{15}], [\Delta_{25}]\rangle\\
\hline

{\mathbb T}^5_{10}& 
\begin{array}{l}
e_1e_2=e_3, e_1e_3=e_4,\\ 
e_2e_4=e_5
\end{array}& \mbox{---} & 
\langle[\Delta_{14}], [\Delta_{23}], [\Delta_{34}]+[\Delta_{15}]\rangle\\
\hline

\end{array}
$$


\subsection{$1$-dimensional  central extensions of ${\mathbb T}^5_{02}$}
Let us use the following notations 
\[\nabla_1=[\Delta_{14}], \nabla_2=[\Delta_{15}], \nabla_3=[\Delta_{23}], \nabla_4=[\Delta_{25}],  \nabla_5=[\Delta_{35}],   
\nabla_6=[\Delta_{24}],\nabla_7= [\Delta_{45}].\]
The automorphism group of ${\mathbb T}^5_{02}$ 
consists of invertible  matrices of the form 
$$\phi=
\begin{pmatrix}
x& 0& 0& 0& 0\\
f& y& 0& 0& 0\\
u& v& xy& 0& 0\\
h& r& xv& x^2y& l\\
t& g& 0& 0& z
\end{pmatrix}.
$$
Notice also that $\det\phi=x^4y^3z\ne 0$.
Since
$$\phi^T\begin{pmatrix}
0& 0& 0& \alpha_1& \alpha_2\\
0& 0& \alpha_3&  \alpha_6 & \alpha_4\\
0& -\alpha_3& 0& 0 &  \alpha_5\\
-\alpha_1& -\alpha_6 & 0 & 0& \alpha_7\\
-\alpha_2& -\alpha_4& -\alpha_5 & -\alpha_7 & 0
\end{pmatrix} \phi=
\begin{pmatrix}
0&\alpha^*&\alpha^{**}&\alpha_1^*&\alpha_2^*\\
-\alpha^*&0&\alpha_3^*&\alpha_6^*&\alpha_4^*\\
-\alpha^{**}&-\alpha_3^*&0&0&\alpha_5^*\\
-\alpha_1^*&-\alpha_6^*&0&0&\alpha_7^*\\
-\alpha_2^*&-\alpha_4^*&-\alpha_5^*&-\alpha_7^*& 0
\end{pmatrix},
$$
where
$$
\begin{array}{rcl} 


\alpha_1^*&=&\alpha_1x^3y+\alpha_6fx^2y-\alpha_7tx^2y,\\

\alpha_2^*&=&\alpha_1 xl+\alpha_2xz+\alpha_4 fz+\alpha_5uz+\alpha_6fl+\alpha_7(hz-tl),\\

\alpha_3^*&=&\alpha_3xy^2-\alpha_5gxy+\alpha_6vxy-\alpha_7gvx,\\

\alpha_4^*&=&\alpha_4 yz+\alpha_5vz+\alpha_6yl+\alpha_7(rz-gl),\\

\alpha_5^* &=&\alpha_5xyz+\alpha_7vxz, \\

\alpha_6^* &=&\alpha_6x^2y^2-\alpha_7gx^2y,\\

\alpha_7^*&=&\alpha_7x^2yz,
\end{array}
$$
we have that the action of $\operatorname{Aut}({\mathbb T}^5_{02})$ on the subspace
$\Big\langle \sum\limits_{i=1}^7 \alpha_i \nabla_i \Big\rangle$
is given by 
$\Big \langle \sum\limits_{i=1}^7 \alpha_i^*  \nabla_i \Big\rangle.$

We are interested in the cocycles with $(\af_6,\af_7)\ne(0,0)$. So, we have the following cases:

\begin{enumerate}
    \item $\alpha_7 \neq 0,$ then choosing
    $g=\frac{\alpha_6 y}{\alpha_7},$
    $v=-\frac{\alpha_5 y}{\alpha_7},$ 
    $l=0,$
    $r=-\frac{\alpha_4y+\alpha_5v}{\alpha_7},$
    $h=-\frac{\alpha_2 x+\alpha_4 f+\alpha_5u}{\alpha_7},$
    $t=\frac{\alpha_1 x +\alpha_6f}{\alpha_7}$ we have the representative $\langle\alpha_3^* \nabla_3+ \alpha_7^* \nabla_7\rangle.$ Here we have two orbits:
    \begin{enumerate}
        \item If $\alpha_3^*\neq 0,$ i.e. $\alpha_3-\frac{\alpha_5\alpha_6}{\alpha_7}\ne 0$, then choosing $y=z=1$, $x=\frac 1{\alpha_7}(\alpha_3-\frac{\alpha_5\alpha_6}{\alpha_7})$ we have the representative $\langle\nabla_3+\nabla_7\rangle.$
        \item If $\alpha_3^*=0,$ i.e. $\alpha_3-\frac{\alpha_5\alpha_6}{\alpha_7}=0$, then we have the representative $\langle\nabla_7\rangle.$
\end{enumerate}

    \item $\alpha_7=0, \alpha_6\neq0$, then choosing
    $l=-\frac{(\alpha_4 y+\alpha_5 v)z}{\alpha_6 y}$,
    $v=\frac{\alpha_5 g-\alpha_3 y}{\alpha_6},$
    $f=-\frac{\alpha_1 x}{\alpha_6}$,
    we have the representative
    $\langle\alpha_2^* \nabla_2+\alpha_5 xyz\nabla_5+\alpha_6x^2y^2 \nabla_6\rangle.$
    Here we have three orbits:
    \begin{enumerate}
        \item If $\alpha_5\neq0,$ then choosing 
        $x=\frac{\alpha_5}{\alpha_6}$, $y=z=1$, $u=-\frac{(\alpha_1l+\alpha_2)x+(\alpha_4+\alpha_6l)f}{\alpha_5}$, 
        we have the representative $\langle\nabla_5+\nabla_6\rangle.$
        
        \item If $\alpha_5=0, \alpha_2^*\ne 0,$ i.e. $\alpha_2-\frac{\alpha_1\alpha_4}{\alpha_6}\ne 0$, then choosing 
        $y=1$, $x=\frac 1{\alpha_6}(\alpha_2-\frac{\alpha_1\alpha_4}{\alpha_6})$,
        we have the representative $\langle\nabla_2+\nabla_6\rangle.$
        
        \item If $\alpha_5=0, \alpha_2^*=0,$ i.e. $\alpha_2-\frac{\alpha_1\alpha_4}{\alpha_6}=0$, then we have the representative $\langle\nabla_6\rangle.$
        
    \end{enumerate}
  
    \end{enumerate}

Note that the representative $\nabla_6$ gives a Tortkara algebra with $2$-dimensional annihilator, which was found in the previous section. A straightforward computation shows that the rest of the representatives belong to different orbits. They give the following pairwise non-isomorphic Tortkara algebras:
$$
\begin{array}{llllllll} 
{\mathbb T}_{04}^6 & : &  
e_1e_2=e_3, & e_1e_3=e_4 , &  e_1e_5=e_6, & e_2e_4=e_6;  \\
{\mathbb T}_{05}^6 & : &  
e_1e_2=e_3, & e_1e_3=e_4 , & e_2e_3=e_6, & e_4e_5=e_6;  \\
{\mathbb T}_{06}^6 & : &  
e_1e_2=e_3, & e_1e_3=e_4, & e_2e_4=e_6, & e_3e_5=e_6;   \\
{\mathbb T}_{07}^6 & : &  
e_1e_2=e_3, & e_1e_3=e_4, & e_4e_5=e_6.   \\
\end{array}$$

\subsection{$1$-dimensional  central extensions of ${\mathbb T}^5_{04}$}
Let us use the following notations 

\[ \nabla_1=[\Delta_{14}],  \ \nabla_2=[\Delta_{25}], 
\ \nabla_3=[\Delta_{15}]+[\Delta_{24}],\ \nabla_4=[\Delta_{15}]. \]
The automorphism group of $\mathbb{T}_{04}^5$ consists of invertible matrices of the form
$$\phi=
\begin{pmatrix}
x& y& 0& 0& 0\\
v& z& 0& 0& 0\\
u& h& xz-yv& 0& 0\\
l& r& xh-yu& x(xz-yv)& y(xz-yv)\\
t& g& vh-zu& v(xz-yv)& z(xz-yv)
\end{pmatrix}.$$
Notice that $\det\phi=(xz-yv)^5\ne 0$.
Since
$$\phi^T\begin{pmatrix}
0&0&0&\alpha_1&\alpha_3+\alpha_4\\
0&0&0&\alpha_3&\alpha_2\\
0&0&0&0&0\\
-\alpha_1&-\alpha_3&0&0&0\\
-\alpha_3-\alpha_4&-\alpha_2&0&0&0
\end{pmatrix}\phi=
\begin{pmatrix}
0&\beta_1&\beta_2&\alpha_1^*&\alpha_3^*+\alpha_4^*\\
-\beta_1&0&\beta_3&\alpha_3^*&\alpha_2^*\\
-\beta_2&-\beta_3&0&0&0\\
-\alpha_1^*&-\alpha_3^*&0&0&0\\
-\alpha_3^*-\alpha_4^*&-\alpha_2^*&0&0&0
\end{pmatrix},$$
where
$$\begin{array}{rcl}



\alpha_1^*&=&(xz-yv)(\af_1 x^2+\af_2 v^2+(2\af_3+\af_4)xv),\\

\alpha_2^* &=&(xz-yv)(\af_1 y^2+\af_2 z^2+(2\af_3+\af_4)yz),\\



\alpha_3^*&=&(xz-yv)(\af_1 xy+\af_2 vz+(\af_3+\af_4)yv+\af_3 xz),\\

\alpha_4^*&=& \alpha_4(xz-yv)^2,

\end{array}
$$
  we have that the action of $\operatorname{Aut}({\mathbb T}^5_{04})$ on the subspace
$\Big\langle \sum\limits_{i=1}^4 \alpha_i \nabla_i \Big\rangle$
is given by 
$\Big \langle \sum\limits_{i=1}^4 \alpha_i^*  \nabla_i\Big\rangle.$

We are only interested in cocycles with $\alpha_4 \neq 0.$ Note that this condition is invariant under automorphisms. We have the following cases:
\begin{enumerate}
    \item $\af_1\ne 0$.\label{af_1-ne-0}
    \begin{enumerate}
        \item $(2\af_3+\af_4)^2-4\af_1\af_2\ne 0$.\label{(2.af_3+af_4)^2-4.af_1.af_2-ne-0} Then for all fixed non-zero $v$ and $z$ the equations
        \begin{align*}
            \af_1 x^2+\af_2 v^2+(2\af_3+\af_4)xv&=0,\\
            \af_1 y^2+\af_2 z^2+(2\af_3+\af_4)yz&=0
        \end{align*}
        have roots $x_1\ne x_2$ and $y_1\ne y_2$. If $x_1z=y_1v$, then $x_2z\ne y_1v$. Thus, we may choose $x,y,z,v$ such that $xz-yv\ne 0$ and $\af_1^*=\af_2^*=0$. Note that $\af_4^*\ne 0$, so applying an automorphism to the cocycle $\af_3^*\nabla_3+\af_4^*\nabla_4$, we may also make $\af_4^*=1$. This gives the representative of the orbit of the form $\xi_{\alpha}=\langle\alpha \nabla_3+\nabla_4\rangle.$ Given a pair $\af,\beta$, it is easy to see that there exists an automorphism $\phi$, such that $\phi(\xi_{\alpha})=\xi_{\beta}$, if and only if $\beta=-\alpha-1.$
        
        \item $(2\af_3+\af_4)^2-4\af_1\af_2=0$. Then choosing $y=-\sqrt{\frac{\af_2}{\af_1}}z$ and $x\ne-\sqrt{\frac{\af_2}{\af_1}}v$, we get the representative $\langle\af_1^*\nabla_1+\af_3^*\nabla_3+\af_4^*\nabla_4\rangle$. A routine calculation shows that the condition $(2\af_3+\af_4)^2-4\af_1\af_2=0$ is invariant under the action of an automorphism, so we have $\af_4^*=-2\af_3^*$. Since $\af_4^*\ne 0$, then there is an automorphism which sends $\af_1^*\nabla_1+\af_3^*\nabla_3-2\af_3^*\nabla_4$ to $\af\nabla_1+\nabla_3-2\nabla_4$ for some $\af\ne 0$. Then choosing $x=\frac 1{\sqrt\af}$, $y=0$, $z=\sqrt{\af}$, we have the representative $\langle\nabla_1+\nabla_3-2\nabla_4\rangle$.
    \end{enumerate}
    
    \item $\af_1=0$, $\af_2\ne 0$. Then choosing $x=z=0$ and $y=v=1$, we get $\alpha_1^*=-\af_2\ne 0$, so we are in \cref{af_1-ne-0}. 
    
    \item $\af_1=0$, $\af_2=0$. This case has already appeared in \cref{(2.af_3+af_4)^2-4.af_1.af_2-ne-0}.
\end{enumerate}
It is easily checked that the orbit of $\langle\nabla_1+\nabla_3-2\nabla_4\rangle$ does not belong to the family of orbits of $\langle\alpha \nabla_3+\nabla_4\rangle$. Thus, we have the following algebras:
$$
\begin{array}{llllllll} 
{\mathbb T}_{08}^6 & : &  
e_1e_2=e_3, & e_1e_3=e_4, & e_1e_4=e_6, & e_1e_5=-e_6,&  e_2e_3=e_5, & e_2e_4=e_6;\\
{\mathbb T}_{09}^6(\alpha) & : &  
e_1e_2=e_3, & e_1e_3=e_4, & e_1e_5=(\alpha+1) e_6,&  e_2e_3=e_5, & e_2e_4=\alpha e_6,
\end{array}
$$
where ${\mathbb T}_{08}^6\not\cong {\mathbb T}_{09}^6(\alpha)$, and ${\mathbb T}_{09}^6(\alpha)\cong {\mathbb T}_{09}^6(\beta)$ if and only if $\beta=-\alpha-1$.


\subsection{$1$-dimensional  central extensions of ${\mathbb T}^5_{06}$}
Let us use the following notations 
\[\nabla_1=[\Delta_{13}], \nabla_2= [\Delta_{14}],  \nabla_3=[\Delta_{24}],  \nabla_4=[\Delta_{25}],  \nabla_5=[\Delta_{34}],  \nabla_6=[\Delta_{15}],  \nabla_7=[\Delta_{45}].\]
The automorphism group of ${\mathbb T}^5_{06}$\ consists of invertible
matrices of the form 
$$\phi=
\begin{pmatrix}
x& p& 0& 0& 0\\
0& y& 0& 0& 0\\
z& t& xy& -py& 0\\
q& r& 0& y^2& 0\\
s& h& xr-yz-pq& f& xy^2
\end{pmatrix}.$$
Notice that $\det\phi=x^3y^6\ne 0$.
Since 
$$\phi^T\begin{pmatrix}
0           &   0           &   \alpha_1    &\alpha_2   &\alpha_6\\
0           &   0           &   0           &\alpha_3   &\alpha_4\\
-\alpha_1   &   0           &   0           &\alpha_5   &0\\
-\alpha_2   &   -\alpha_3   &   -\alpha_5   &0          &\alpha_7\\
-\alpha_6   &   -\alpha_4   &   0           &-\alpha_7  &0
\end{pmatrix}\phi=
\begin{pmatrix}
0                       &   \alpha^*     &\alpha_1^*     &\alpha_2^*+\alpha^{**} &\alpha_6^*\\
-\alpha^*                &   0           &\alpha^{**}        &\alpha_3^*&\alpha_4^*\\
-\alpha_1^*             &   -\alpha^*    &0&\alpha_5^*   &0\\
-\alpha_2^* -\alpha^{**}    &   -\alpha_3^* &-\alpha_5^*    &0&\alpha_7^*\\
-\alpha_6^*             &   -\alpha_4^* &0              &-\alpha_7^*&0
\end{pmatrix},
$$
where
$$\begin{array}{rcl}

\alpha_1^* &=& xy(\af_1x-\af_5q) +(\af_6x+\af_7q)(rx-yz-pq),\\


\alpha_2^*&=&-2\af_1pxy+\alpha_2xy^2+\af_4y(pq-rx+yz)+\af_5y(pq+rx+yz)\\
&&+\af_6(p(pq-rx+yz)+fx)+\af_7(r(pq-rx+yz)-sy^2+fq),\\


\alpha_3^*&=& 
y(-\alpha_1p^2+\alpha_2py+\alpha_3y^2+\alpha_4f)+\alpha_5y(pr+yt)+\alpha_6fp+\alpha_7(fr-y^2h),\\



\alpha_4^*&=&xy^2(\alpha_4y+\alpha_6p+\alpha_7r),\\

\alpha_5^*&=& y^2(\af_5xy+ \af_7(pq-rx+yz)),\\

\alpha_6^*&=&xy^2(\alpha_6x+\alpha_7q), \\

\alpha_7^*&=& \alpha_7y^4x,
\end{array}
$$
we have that the action of $\operatorname{Aut}({\mathbb T}^5_{06})$ on the subspace
$\Big\langle \sum\limits_{i=1}^7 \alpha_i \nabla_i \Big\rangle$
is given by 
$\Big \langle \sum\limits_{i=1}^7 \alpha_i^*  \nabla_i\Big\rangle.$
 We are only interested in the cocycles with $(\af_6,\af_7)\ne(0,0)$. There are the following cases:

\begin{enumerate}
   \item  $\alpha_7\neq 0.$ 
    Then choosing 
    $q=-\frac{\alpha_6}{\alpha_7}x$, 
    $z= -\frac{\alpha_5xy+ \af_7(pq-rx)}{\alpha_7y}$, 
    $r=-\frac{\alpha_4y+\alpha_6p}{\alpha_7}$, 
     and $h,s$ such that $\alpha_2^*=\alpha_3^*=0,$
     we have $\af_1^*=\frac{\alpha_1\alpha_7+\alpha_5\alpha_6}{\af_7}x^2y$. So, there are 2 subcases:
    
    \begin{enumerate}
        \item $\alpha_1\alpha_7+\alpha_5\alpha_6\neq 0$. Then choose $x=\sqrt[7]{\frac{\af_7^5}{(\alpha_1\alpha_7+\alpha_5\alpha_6)^4}}$ and $y=\sqrt[7]{\frac{\alpha_1\alpha_7+\alpha_5\alpha_6}{\af_7^3}}$. This gives the representative $\langle \nabla_1 +\nabla_7 \rangle$.
        \item $\alpha_1\alpha_7+\alpha_5\alpha_6=0$. Then choose $x=\frac 1{\af_7}$, $y=1$ to get the representative $\langle \nabla_7 \rangle.$
    \end{enumerate}
\item $\alpha_7= 0$, $\af_6\ne 0$.
 Then choosing $r$ and $f$, we may make $\af_1^*=\af_2^*=0$. Moreover, if $p=-\frac{\af_4}{\af_6}y$, then $\af_4^*=0$. For this value of $p$ we have $\af_3^*=y^3(\af_3-\frac{\af_1\af_4^2}{\af_6^2}-\frac{\af_2\af_4}{\af_6})+\af_5y^2(t-\frac{\af_4}{\af_6}r)$. So, we have two subcases:
    \begin{enumerate}
        \item $\af_5\ne 0$. Then choosing the appropriate value of $t$, we may make $\af_3^*=0$. So, $x=\frac{\sqrt{\af_5}}{\sqrt[4]{\af_6^3}}$ and $y=\frac{\sqrt[4]{\af_6}}{\sqrt{\af_5}}$ give the representative $\langle  \nabla_5+ \nabla_6 \rangle$.
        \item $\af_5=0$. We have two subcases:
        \begin{enumerate}
            \item $\af_3\af_6^2-\af_1\af_4^2-\af_2\af_4\af_6\ne 0$. Then choosing $y=\sqrt[3]{\frac{\af_6^2}{\af_3\af_6^2-\af_1\af_4^2-\af_2\af_4\af_6}}$ and $x=\frac{\sqrt[3]{\af_3\af_6^2-\af_1\af_4^2-\af_2\af_4\af_6}}{\sqrt[6]{\af_6}}$ we get the representative $\langle \nabla_3+\nabla_6 \rangle$.
            \item $\af_3\af_6^2-\af_1\af_4^2-\af_2\af_4\af_6=0$. Then choosing $x=1$, $y=\frac 1{\sqrt{\af_6}}$, we get the representative $\langle \nabla_6 \rangle$.
        \end{enumerate}
    \end{enumerate}



\end{enumerate}
A straightforward verification shows that the found orbits are distinct. Thus, we have new non-isomorphic Tortkara algebras:
$$
\begin{array}{lllllllll} 
{\mathbb T}_{10}^6 & : &  
e_1e_2=e_3, & e_1e_3=e_6,  & e_1e_4=e_5, & e_2e_3=e_5, & e_4e_5=e_6; \\

{\mathbb T}_{11}^6 & :& 
e_1e_2=e_3, & e_1e_4=e_5, & e_1e_5=e_6, &  e_2e_3=e_5;     \\

{\mathbb T}_{12}^6 & :& 
e_1e_2=e_3, & e_1e_4=e_5, & e_1e_5=e_6, &  e_2e_3=e_5, &  e_2e_4=e_6;  \\

{\mathbb T}_{13}^6 & :& 
e_1e_2=e_3, & e_1e_4=e_5, & e_1e_5=e_6, &  e_2e_3=e_5, & e_3e_4=e_6; \\

{\mathbb T}_{14}^6 & : & 
e_1e_2=e_3, & e_1e_4=e_5, & e_2e_3=e_5, & e_4e_5=e_6.

\end{array}$$


\subsection{$1$-dimensional  central extensions of ${\mathbb T}^5_{08}$}
Let us use the following notations 
\[\nabla_1=[\Delta_{15}], \nabla_2= [\Delta_{23}],  \nabla_3=[\Delta_{24}], \nabla_4=[\Delta_{25}].\]
The automorphism group of ${\mathbb T}^5_{08}$\ consists of invertible
matrices of the form 
$$\phi=
\begin{pmatrix}
x& 0& 0& 0& 0\\
z& y& 0& 0& 0\\
t& p& xy& 0& 0\\
q& r& xp& x^2y& 0\\
h&s&xr&x^2p&x^3y\end{pmatrix}.$$
Notice that $\det\phi=x^7y^4\ne 0$. Since 

$$\phi^T\begin{pmatrix}
0&0&0&0&\alpha_1\\
0&0&\alpha_2&\alpha_3&\alpha_4\\
0&-\alpha_2&0&0&0\\
0&-\alpha_3&0&0&0\\
-\alpha_1&-\alpha_4&0&0&0
\end{pmatrix}\phi=
\begin{pmatrix}
0&\beta_1&\beta_2&\beta_3&\alpha_1^*\\
-\beta_1&0&\alpha_2^*&\alpha_3^*&\alpha_4^*\\
-\beta_2&-\alpha_2^*&0&0&0\\
-\beta_3&-\alpha_3^*&0&0&0\\
-\alpha_1^*&\alpha_4^*&0&0&0
\end{pmatrix},
$$
where
$$\begin{array}{rcl}
\alpha_1^* &=&x^3y(\alpha_1x+\alpha_4z), \\

\alpha_2^* &=&xy(\alpha_2y+\alpha_3p+\alpha_4r), \\

\alpha_3^* &=& x^2y(\alpha_3y+\alpha_4p), \\

\alpha_4^* &=& \alpha_4 x^3y^2,
\end{array}$$
we have that the action of $\operatorname{Aut}({\mathbb T}^5_{08})$ on the subspace
$\Big\langle \sum\limits_{i=1}^4 \alpha_i \nabla_i \Big\rangle$
is given by 
$\Big\langle \sum\limits_{i=1}^4 \alpha_i^*  \nabla_i \Big\rangle.$
 We are interested in cocycles with $(\af_3,\af_4)\ne(0,0)$.
 So, we have the following cases:
 \begin{enumerate}
     \item $\alpha_4\neq0.$ Then choosing $p=-\frac{\alpha_3x}{\alpha_4}, r=-\frac{\alpha_2 y+\alpha_3p}{\alpha_4}$ and
     $z=-\frac{\alpha_1x}{\alpha_4}$ we have the representative $\langle\nabla_4\rangle.$

     \item $\alpha_4=0, \alpha_3 \neq 0$. We have two subcases:
     \begin{enumerate}
         \item $\alpha_1\neq 0$. Then choosing 
     $p=-\frac{\alpha_2y}{\alpha_3}$ and $y=\frac{\alpha_1 x^2}{\alpha_3}$ we have the representative $\langle\nabla_1+\nabla_3\rangle.$
     
        \item $\alpha_1=0$. Then choosing $p=-\frac{\alpha_2y}{\alpha_3}$ we have the representative $\langle\nabla_3\rangle.$
     \end{enumerate}
 \end{enumerate}
Note that the representative $\langle\nabla_3\rangle$ gives a
$6$-dimensional Tortkara algebra with $2$-dimensional annihilator.
It is a $2$-dimensional central extension of a $4$-dimensional Tortkara algebra, and it was found in the previous section.
The rest of the representatives give new non-isomorphic Tortkara algebras:
$$
\begin{array}{lllllllll} 
{\mathbb T}_{15}^6 & : &  
e_1e_2=e_3, & e_1e_3=e_4, & e_1e_4=e_5, & e_1e_5=e_6, & e_2e_4=e_6; \\

{\mathbb T}_{16}^6 & : &  
e_1e_2=e_3, & e_1e_3=e_4, & e_1e_4=e_5, & e_2e_5=e_6.

\end{array}$$


\subsection{$1$-dimensional  central extensions of ${\mathbb T}^5_{09}$}
Let us use the following notations 
\[\nabla_1=[\Delta_{14}], \nabla_2= [\Delta_{24}]+[\Delta_{15}],  \nabla_3=[\Delta_{15}], \nabla_4=[\Delta_{25}].\]
The automorphism group of ${\mathbb T}^5_{09}$\ consists of invertible
matrices of the form 
$$\phi=
\begin{pmatrix}
x& 0& 0& 0& 0\\
y& x^2& 0& 0& 0\\
z& t& x^3& 0& 0\\
p& q& xt& x^4& 0\\
r&s&-x^2z+xq+yt&x^3y+x^2t&x^5\end{pmatrix}.$$
Notice that $\det\phi=x^{15}\ne 0$. Since 
$$\phi^T\begin{pmatrix}
0&0&0&\alpha_1&\alpha_2+\alpha_3\\
0&0&0&\alpha_2&\alpha_4\\
0&0&0&0&0\\
-\alpha_1&-\alpha_2&0&0&0\\
-\alpha_2-\af_3&-\alpha_4&0&0&0
\end{pmatrix}\phi =
\begin{pmatrix}
0                           &\alpha^{*}     &\alpha^{**}    &\alpha^{***}+\alpha_1^*&\alpha_2^*+\alpha_3^*\\
-\alpha^{*}                 &0              &\alpha^{***}   &\alpha_2^*&\alpha_4^*\\
-\alpha^{**}                &-\alpha^{***}  &0              &0&0\\
-\alpha^{***}-\alpha_1^*    &\alpha_2^*     &0              &0&0\\
-\alpha_2^*-\alpha_3^*      &-\alpha_4^*    &0              &0&0
\end{pmatrix},$$
where
$$
\begin{array}{rcl}
\alpha_1^*&=&x^3(\af_1x^2 + 2\af_2xy + \af_3(xy+t) + \af_4(y^2 + xz - q)),\\

\alpha_2^* &=&x^4(\alpha_2x^2+\alpha_4(xy+t)),\\

\alpha_3^* &=& x^4(\alpha_3x^2 -\alpha_4t),\\

\alpha_4^*&=& \alpha_4 x^7,
\end{array}
$$
we have that the action of $\operatorname{Aut}({\mathbb T}^5_{09})$ on the subspace
$\Big\langle \sum\limits_{i=1}^4 \alpha_i \nabla_i \Big\rangle$
is given by 
$\Big\langle \sum\limits_{i=1}^4 \alpha_i^*  \nabla_i \Big\rangle.$
We are only interested in cocycles with $(\af_3,\af_4)\ne(0,0)$. So, we have the following cases:

\begin{enumerate}
    \item $\alpha_4\neq0$. Then choosing 
    $t=\frac{\alpha_3 x^2}{\alpha_4}, y=-\frac{\alpha_2x^2+\alpha_4t}{\alpha_4x}$ and 
    $z=-\frac{\af_1x^2 + 2\af_2xy + \af_3(xy+t) + \af_4(y^2 - q)}{\alpha_4x}$ we have the representative $\langle\nabla_4\rangle.$

\item $\alpha_4=0, \alpha_3 \neq 0$. Then choosing 
   $t=-\frac{\af_1x^2 + 2\af_2xy + \af_3xy}{\alpha_3}$ we have the family of representatives $\langle\alpha\nabla_2+\nabla_3\rangle.$ Observe that there exists an automorphism $\phi$, such that $\phi\langle\alpha\nabla_2+\nabla_3\rangle=\langle\beta\nabla_2+\nabla_3\rangle$ if and only if $\af=\beta$.
\end{enumerate}

Thus, we have the following family of pairwise non-isomorphic Tortkara algebras:
$$
\begin{array}{llllllllll} 

\mathbb{T}_{17}^6 &:&
e_1e_2=e_3, & e_1e_3=e_4, & e_1e_4=e_5, & e_2e_3=e_5, & e_2e_5=e_6;  \\

\mathbb{T}_{18}^6(\alpha) &:&
e_1e_2=e_3, & e_1e_3=e_4, & e_1e_4=e_5, &e_1e_5=(\alpha+1)e_6, & e_2e_3=e_5, & e_2e_4=\alpha e_6.  
\end{array}
$$


\subsection{$1$-dimensional  central extensions of ${\mathbb T}^5_{10}$}
Let us use the following notations 
\[\nabla_1=[\Delta_{14}], \nabla_2= [\Delta_{23}],  \nabla_3=[\Delta_{15}]+[\Delta_{34}].\]
The automorphism group of ${\mathbb T}^5_{10}$\ consists of invertible
matrices of the form 
$$\phi=
\begin{pmatrix}
x& 0& 0& 0& 0\\
0& y& 0& 0& 0\\
z& 0& xy& 0& 0\\
p&q& 0& x^2y& 0\\
r& h& -yp& 0& y^2x^2
\end{pmatrix}.$$
Notice that $\det\phi=x^6y^5\ne 0$.
Since
$$\phi^T
\begin{pmatrix}
0& 0& 0& \alpha_1& \alpha_3\\
0& 0& \alpha_2& 0& 0\\
0& -\alpha_2& 0& \alpha_3& 0\\
-\alpha_1& 0& -\alpha_3& 0& 0\\
-\alpha_3& 0& 0& 0& 0
\end{pmatrix} \phi=
\begin{pmatrix}
0&\alpha^*& \alpha^{**}& \alpha_1^*&\alpha_3^*\\
-\alpha^*& 0& \alpha_2^*&0&0\\
-\alpha^{**}&-\alpha_2^*&0&\alpha_3^*&0\\
-\alpha_1^*&0&-\alpha_3^*&0&0\\
-\alpha_3^*&0&0&0&0
\end{pmatrix},$$
where
$$\begin{array}{rcl}
\alpha_1^* &=& x^2y(\alpha_1x+ \alpha_3z),\\
\alpha_2^* &=& xy(\alpha_2y-\alpha_3q),\\
\alpha_3^* &=& \alpha_3 x^3y^2,
\end{array}$$
we have that the action of $\operatorname{Aut}({\mathbb T}^5_{10})$ on the subspace
$\Big\langle \sum\limits_{i=1}^3 \alpha_i \nabla_i \Big\rangle$
is given by 
$\Big\langle \sum\limits_{i=1}^3 \alpha_i^*\nabla_i \Big\rangle.$

Every element with $\alpha_3=0$ gives an algebra with $2$-dimensional annihilator, which was found in the previous section.
If $\alpha_3\neq0,$ then choosing $x=\frac{1}{\sqrt[3]{\alpha_3}}, y=1, z=-\frac{\alpha_1x}{\alpha_3}, q=\frac{\alpha_2}{\alpha_3}$ we have the representative $\langle\nabla_3\rangle.$

Thus, we obtain the following new Tortkara algebra:
$$\begin{array}{lllll llll}
{\mathbb T}_{19}^6&:& 
e_1e_2=e_3,& e_1e_3=e_4,& e_1e_5=e_6, & e_2e_4=e_5, & e_3e_4=e_6.  
\end{array}$$



\section{Main result}
The algebraic classification of all $6$-dimensional nilpotent Malcev algebras was obtained in \cite{hac18}. As it was proved in \cite{hac18}, every $6$-dimensional nilpotent Malcev algebra is metabelian (i.e. $(xy)(zt)=0$). In particular, every $6$-dimensional nilpotent Lie algebra is metabelian, and hence it is a Tortkara algebra.
By some easy verification of all $6$-dimensional nilpotent non-Lie Malcev algebras, one can see that all these algebras are also Tortkara.
Now we have

\begin{lemma}
 Let $\mathbb T$ be a $6$-dimensional nilpotent Malcev-Tortkara algebra over $\mathbb C$ with non-trivial product. Then $\mathbb T$ is isomorphic to one of the nilpotent Malcev algebras listed in \cite{hac18}.
\end{lemma}

Now we are ready to formulate the main result of our paper:

\begin{theorem}
Let $\mathbb T$ be a $6$-dimensional nilpotent non-Malcev Tortkara algebra over $\mathbb C$. Then $\mathbb T$ is isomorphic to one of the following algebras:

\begin{enumerate}[$\bullet$]

\item ${\mathbb T}_{00}^6 \ : \  e_1e_2=e_3, \  e_1e_3=e_4, \ e_2e_4=e_5;$

\item ${\mathbb T}_{01}^6 \ : \  
e_1e_2=e_3,\ e_1e_3=e_4,\ e_1e_4=e_5,\ e_2e_3=e_5,\ e_2e_4=e_6;$

\item ${\mathbb T}_{02}^6 \ : \ 
e_1e_2=e_3,\ e_1e_3=e_4,\ e_2e_3=e_5,\ e_2e_4=e_6;$

\item ${\mathbb T}_{03}^6 \ : \    
e_1e_2=e_3, \ e_1e_3=e_4, \ e_1e_4=e_5, \ e_2e_4=e_6;$

\item ${\mathbb T}_{04}^6 \ : \    
e_1e_2=e_3, \ e_1e_3=e_4, \ e_1e_5=e_6, \ e_2e_4=e_6;$ 

\item ${\mathbb T}_{05}^6 \ : \
e_1e_2=e_3, \ e_1e_3=e_4 , \ e_2e_3=e_6, \ e_4e_5=e_6;$
\item ${\mathbb T}_{06}^6 \ : \
e_1e_2=e_3, \ e_1e_3=e_4, \ e_2e_4=e_6, \ e_3e_5=e_6  ;$

\item ${\mathbb T}_{07}^6 \ : \
e_1e_2=e_3, \ e_1e_3=e_4, \  e_4e_5=e_6 ;$

\item ${\mathbb T}_{08}^6 \ : \  
e_1e_2=e_3, \ e_1e_3=e_4, \ e_1e_4=e_6, \ e_1e_5=-e_6, \ e_2e_3=e_5, \ e_2e_4=e_6;$
\item ${\mathbb T}_{09}^6(\alpha) \ : \  
e_1e_2=e_3, \ e_1e_3=e_4, \ e_1e_5=(\alpha+1) e_6, \ e_2e_3=e_5, \ e_2e_4=\alpha e_6;$

\item ${\mathbb T}_{10}^6 \ : \  
e_1e_2=e_3, \ e_1e_3=e_6,  \ e_1e_4=e_5, \ e_2e_3=e_5, \ e_4e_5=e_6;$

\item ${\mathbb T}_{11}^6 \ : \ 
e_1e_2=e_3, \ e_1e_4=e_5, \ e_1e_5=e_6, \  e_2e_3=e_5;$

\item ${\mathbb T}_{12}^6 \ : \ 
e_1e_2=e_3, \ e_1e_4=e_5, \ e_1e_5=e_6, \  e_2e_3=e_5, \  e_2e_4=e_6;$

\item ${\mathbb T}_{13}^6 \ : \ 
e_1e_2=e_3, \ e_1e_4=e_5, \ e_1e_5=e_6, \  e_2e_3=e_5, \ e_3e_4=e_6;$

\item ${\mathbb T}_{14}^6 \ : \ 
e_1e_2=e_3, \ e_1e_4=e_5, \ e_2e_3=e_5, \ e_4e_5=e_6;$

\item ${\mathbb T}_{15}^6 \ : \ 
e_1e_2=e_3, \ e_1e_3=e_4, \ e_1e_4=e_5, \ e_1e_5=e_6, \ e_2e_4=e_6 ;$

\item ${\mathbb T}_{16}^6 \ : \  
e_1e_2=e_3, \ e_1e_3=e_4, \ e_1e_4=e_5, \ e_2e_5=e_6;$

\item $\mathbb{T}_{17}^6 \ : \
e_1e_2=e_3, \ e_1e_3=e_4, \ e_1e_4=e_5, \ e_2e_3=e_5, \ e_2e_5=e_6 ;$

\item $\mathbb{T}_{18}^6(\alpha) \ : \
e_1e_2=e_3, \ e_1e_3=e_4, \ e_1e_4=e_5, \ e_1e_5=(\alpha+1)e_6, \ e_2e_3=e_5, \ e_2e_4=\alpha e_6;$

\item ${\mathbb T}_{19}^6 \ : \ 
e_1e_2=e_3, \ e_1e_3=e_4, \ e_1e_5=e_6, \ e_2e_4=e_5, \ e_3e_4=e_6.$

\end{enumerate}

All listed algebras are non-isomorphic, except  
 ${\mathbb T}^6_{09}(\alpha) \cong {\mathbb T}^6_{09}(-\alpha-1).$

\end{theorem}

\begin{remark}
It was proved that all  $6$-dimensional nilpotent Malcev algebras are metabelian (i.e. $(xy)(zt)=0$) \cite{hac18}.
As follows from the main theorem, 
there is only one $6$-dimensional nilpotent non-metabelian  Tortkara algebra: ${\mathbb T}_{19}^6$. 
\end{remark}

\section{Appendix: The classification of $3$-dimensional Tortkara algebras}

Recall that the complete algebraic, geometric and degeneration classification of all anticommutative 
$3$-dimensional algebras was given in \cite{ikv18}.
Easy corollary gives the algebraic and geometric classification of all $3$-dimensional Tortkara algebras.

\begin{theorem}
Let $\mathfrak{Tort}_3$ be the variety of all $3$-dimensional Tortkara algebra over $\mathbb C$
and $\mathbb T$ be a nonzero algebra from $\mathfrak{Tort}_3,$ then

\begin{enumerate}
    \item  $\mathbb T$ is isomorphic to one of the following algebras:

\begin{enumerate}[$\bullet$]
\item 
$\mathfrak{g}_{1}: e_2e_3 =e_1$

\item 
$ \mathfrak{g}_2 :
e_1e_3 =e_1, \  e_2e_3=e_2 $ 

\item $\mathfrak{g}^{\alpha}_3: 
e_1e_3 =e_1+e_2, \  e_2e_3=\alpha e_2$

\item 
$\mathcal{A}_2: e_1e_2=e_1, \ e_2e_3=e_2$

\item 
$\mathcal{A}_1^{0}: e_1e_2=e_3, \ e_1e_3 =e_1+e_3, \  e_2e_3=\alpha e_2$
\end{enumerate}

\item The variety of $3$-dimensional Tortkara algebras has one irreducible component defined 
by  the rigid algebra $\mathcal{A}_1^{0}.$

\end{enumerate}

\end{theorem}

\end{document}